\documentclass[10pt]{amsart}
\usepackage{mathptmx}
\usepackage{amsmath}
\usepackage{amssymb}
\usepackage{array}
\usepackage{geometry}
\usepackage[bookmarks=true,colorlinks=true, pdfstartview=FitV, linkcolor=black, citecolor=blue, urlcolor=black]{hyperref}

\usepackage{color}
\definecolor{DarkRed}{rgb}{0.55,.00,0.2}
\definecolor{DarkGrey}{rgb}{0.35,.35,0.35}

\theoremstyle{definition}

\theoremstyle{remark}

\numberwithin{equation}{section}

\begin{document}

\title{A new  Weber  type integral equation \\ related to the Weber-Titchmarsh problem }

\author{S. Yakubovich}
\address{Department of Mathematics, Faculty of Sciences,  University of Porto,  Campo Alegre str.,  687; 4169-007 Porto,  Portugal}
\email{ syakubov@fc.up.pt}


\date{\today}
\maketitle

\markboth{\rm \centerline{ S.  Yakubovich}}{}
\markright{\rm \centerline{Weber type integral equation}}

\begin{abstract}   We derive  solvability conditions and closed-form solution for  the Weber type integral equation, related to the familiar Weber-Orr integral transforms and  the old Weber-Titchmarsh problem (posed in {\it Proc. Lond. Math. Soc.}  {\bf 22}(2) (1924),  pp.15, 16), recently solved by the author.    Our method involves properties of the inverse Mellin transform  of integrable functions.  The Mellin-Parseval equality and some integrals, involving the Gauss hypergeometric function  are used.  
\end{abstract}

\vspace{1cm}

Recently, the author gave solvability conditions and closed-form solution for the classical Weber equation \cite{yakweb}

$$ \int_0^\infty   C_\nu (x \xi, a\xi) g(\xi) d\xi = f(x),\eqno(1)$$
where
$f(x)$ is a given function on $[a, \infty),\ a >0,  g(x), x \in \mathbb{R}_+$  should be determined and the kernel 
$$C_\nu(\alpha, \beta)=  J_\nu(\alpha) Y_\nu(\beta) -   Y_\nu(\alpha) J_\nu(\beta)\eqno(2)$$
involves  Bessel functions of the first and second kind  $J_\nu(z), Y_\nu(z),\  \nu \in \mathbb{C}$   \cite{erd}, Vol. II .
It  was solved formally by Titchmarsh in 1924 and posed as an open problem (see  \cite{titweb}, p. 15.)  to describe a class of  complex-valued functions $g(x),\ x \in \mathbb{R}_+$, which  can be expanded in terms of  the following repeated  integral 

$$g(x)=  {x\over  J_\nu^2(ax)+ Y_\nu^2(ax)} \int_a^\infty C_\nu (xt, xa) t \int_0^\infty   C_\nu (t\xi, a\xi) g(\xi) d\xi dt,\ x >0.\eqno(3)$$  
Expansion (3) is related to the familiar Weber-Orr integral expansions of an arbitrary function $f(x)$ as repeated integrals

$$f(x)=  \int_0^\infty  {t \  C_\nu (xt, at)\over  J_\nu^2(at)+ Y_\nu^2(at)}\int_a^\infty   C_\nu (\xi t, at) \xi f(\xi) d\xi dt,\eqno(4)$$  

$$f(x)=   \int_a^\infty C_\nu (xt, xa) t \int_0^\infty    { C_\nu (t \xi, a\xi) \over  J_\nu^2(a\xi)+ Y_\nu^2(a\xi)} \xi f (\xi) d\xi dt,\eqno(5)$$  
which are different from (3).    Our method is  based on the use of the Mellin transform \cite{tit}.  Precisely, the Mellin transform is defined in  $L_{\mu, p}(\mathbb{R}_+),\ 1 < p \le 2$  by the integral  
$$f^*(s)= \int_0^\infty f(x) x^{s-1} dx,\eqno(6)$$
 being convergent  in mean with respect to the norm in $L_q(\mu- i\infty, \mu + i\infty),\   q=p/(p-1)$.   Moreover, the  Parseval equality holds for $f \in L_{\mu, p}(\mathbb{R}_+),\  g \in L_{1-\mu, q}(\mathbb{R}_+)$
$$\int_0^\infty f(x) g(x) dx= {1\over 2\pi i} \int_{\mu- i\infty}^{\mu+i\infty} f^*(s) g^*(1-s) ds.\eqno(7)$$
The inverse Mellin transform is given accordingly
 $$f(x)= {1\over 2\pi i}  \int_{\mu- i\infty}^{\mu+i\infty} f^*(s)  x^{-s} ds,\eqno(8)$$
where the integral converges in mean with respect to the norm  in   $L_{\mu, p}(\mathbb{R}_+)$
$$||f||_{\mu,p} = \left( \int_0^\infty  |f(x)|^p x^{\mu p-1} dx\right)^{1/p}.\eqno(9)$$
In particular, letting $\mu= 1/p$ we get the usual space $L_1(\mathbb{R}_+)$.   We will modify the definition of a  special class of functions related to the Mellin transform (6) and its inversion (8),  which was introduced in \cite{class},\ \cite{yal}. Indeed, we have

{\bf Definition 1}. Denote by ${\mathcal M}^{-1}(L_c)$ the space of
functions $f(x), \ x \in \mathbb{R}_+$,  being representable by inverse
Mellin transform (8) of  integrable functions $F(s) \in L_{1}(c)$ on
the vertical line $c =\{s \in \mathbb{C}:  \mu ={\rm Re s} = c_0\}$.

The  space ${\mathcal M}^{-1}(L_c)$  with  the  usual operations  of
addition   and multiplication by scalar is a linear vector space. If
the norm in ${\mathcal M}^{-1}(L_c)$ is introduced by the formula
$$ \big\vert\big\vert f \big\vert\big\vert_{{\mathcal
M}^{-1}(L_c)}= {1\over 2\pi }\int^{+\infty}_{-\infty} |
F\left(c_0 +it\right)| dt,\eqno(10)$$
then it becomes  a Banach space.

 {\bf Definition 2 }.  Let $\mu\neq 0,\ c_1, c_2 \in \mathbb{R}$ be such that $2 \hbox{sign}\ c_1 + \hbox{sign}\  c_2 \ge 0$. By ${\mathcal M}_{c_1,c_2}^{-1}(L_c)$ we denote the space of functions $f(x), x \in \mathbb{R}_+$, representable in the form (8), where $s^{c_2}e^{\pi c_1|s|} F(s) \in L_1(c)$.

It is a Banach space with the norm
$$ \big\vert\big\vert f \big\vert\big\vert_{{\mathcal
M}_{c_1,c_2}^{-1}(L_c)}= {1\over 2\pi }\int_{c} e^{\pi c_1|s|}
|s^{c_2} F(s) ds|.$$
In particular, letting $c_1=c_2=0$ we get the space ${\mathcal M}^{-1}(L_c)$. Moreover, it is easily seen the inclusion

$${\mathcal M}_{d_1,d_2}^{-1}(L_c) \subseteq {\mathcal
M}_{c_1,c_2}^{-1}(L_c)$$ when $2 \hbox{sign}(d_1- c_1) + \hbox{sign}
(d_2-c_2) \ge 0$.

Using this technique we proved the following 

{\bf Theorem 1\  \cite{yakweb}.} {\it Let $a > 0, \ \nu \in \mathbb{C}, \   0< {\rm Re } \nu < 1/2,  g(x)  \in  {\mathcal M}_{0,1} ^{-1}(L_c)$ with $c =\{s \in \mathbb{C}:  -1< {\rm Re} s < 0\}$.   Then for almost all $x >0$ expansion $(3)$ holds, where the inner and outer integrals are understood in the  improper sense.} 

These results will be applied to solve the so-called Weber type integral equation

$$\int_0^\infty \varphi(\lambda) \left[ J_\nu(x\lambda) Y_{\nu+1} (a\lambda) -   Y_\nu(x\lambda) J_{\nu+1} (a\lambda)\right] d\lambda = f(x),\ x >  a > 0\eqno(11)$$
in the class  $ {\mathcal M}_{0,1} ^{-1}(L_c)$.   The key ingredient will be  also properties for derivative of Bessel functions, namely,  (see \cite{erd}, Vol. II)

$$\left[x^{\mp \nu}  {d\over dx}  x^{\pm \nu} \right] J_\nu(x) = \pm \  J_{\nu \mp1} (x),\eqno(12)$$

$$\left[ x^{\mp \nu}  {d\over dx}  x^{\pm \nu} \right] Y_\nu(x) = \pm \  Y_{\nu\mp 1} (x),\eqno(13)$$ 
and the following integral, which is a direct consequence of relation (2.13.15.4) in \cite{prud}, Vol. 2, namely, 

$$F_\nu(x,s)= \int_0^\infty \lambda^{-s} \left[ J_\nu(x\lambda) Y_{\nu+1} (a\lambda) -   Y_\nu(x\lambda) J_{\nu+1} (a\lambda)\right] d\lambda=   {2^{-s} a^{\nu+1} \over \pi x^{2-s +\nu} }{\cos(\pi\nu)\over \Gamma(s/2)} \Gamma( -\nu-1) \Gamma( 1+\nu- s/2)$$

$$\times {}_2F_1\left( 1- {s\over 2},\  1+\nu- {s\over 2}; \ 2+\nu;\  {a^2\over x^2}\right) -  {2^{-s} x^{\nu+s} \over \pi a^{1 +\nu} }{ \Gamma( \nu+1) \Gamma( - s/2)\over \Gamma(1+\nu + s/2)}$$

$$\times {}_2F_1\left( -\nu -  {s \over 2},\  - {s\over 2}; \  -\nu;\  {a^2\over x^2}\right)-  {2^{-s} a^{\nu+1} \over \pi x^{1 +\nu-s} } \cos \left({\pi s\over 2}\right)  { \Gamma( \nu+1- s/2) \Gamma(1  - s/2)\over \Gamma(2+\nu)}$$

$$\times {}_2F_1\left(1+\nu -  {s \over 2},\  1 - {s\over 2}; \ 2+\nu;\  {a^2\over x^2}\right),\  -1 < {\rm Re} s < 0,\ x >  a,\eqno(14)$$
where $\Gamma(z)$ is Euler's gamma-function and ${}_2F_1(a,b; c; x)$ is Gauss's hypergeometric function \cite{erd}, Vol. 1, having an integral representation as the Euler integral

$$ {}_2F_1(a,b; c; x) = {\Gamma(c)\over \Gamma(b) \Gamma(c-b)} \int_0^1 u^{b-1} (1-u)^{c-b-1} (1- xu)^{-a} du,\eqno(15)$$
for instance, under conditions ${\rm Re} a >0,\ {\rm Re}  c >  {\rm Re}  b > 0\   x \in[ 0,1).$   Moreover, representation (15) gives us the following uniform estimate for the Gauss function, which will be used below

$$\left| {}_2F_1(a,b; c; x) \right| \le \left|  {\Gamma(c)\over \Gamma(b) \Gamma(c-b)} \right| \int_0^1 u^{{\rm Re} b-1} (1-u)^{{\rm Re} (c-b)-1} (1- xu)^{-{\rm Re} a} du$$ 

$$\le   (1- x)^{- {\rm Re} a}\  B\left( {\rm Re} b,\ {\rm Re} (c-b) \right)   \left|  {\Gamma(c)\over \Gamma(b) \Gamma(c-b)} \right|,\eqno(16)$$
where $B(a,b)$ is Euler's beta-function \cite{erd}, Vol. 1.    Let  $-1<  {\rm Re } \nu < -1/2$. Then, denoting by the same letter $C$ various positive constants, which can occur, we obtain 

$$\left| {}_2F_1\left( 1- {s\over 2},\  1+\nu- {s\over 2}; \ 2+\nu;\  {a^2\over x^2}\right)\right| \le C\    x^{2-  {\rm Re} s}  (x^2- a^2)^{{\rm Re} s/2 -1} \left|  {\Gamma(2+\nu)\over \Gamma( (2(1+\nu)- s)/ 2) \Gamma((2+s)/2)} \right|,$$

$$\left| {}_2F_1\left( -\nu -  {s \over 2},\  - {s\over 2}; \  -\nu;\  {a^2\over x^2}\right)\right| \le C\    x^{-2\nu -  {\rm Re} s}  (x^2- a^2)^{\nu+ {\rm Re} s/2} \left|  {\Gamma(-\nu)\over \Gamma( - s/2) \Gamma( -\nu+ s/ 2)} \right|,$$

$$\left|{}_2F_1\left(1+\nu -  {s \over 2},\  1 - {s\over 2}; \ 2+\nu;\  {a^2\over x^2}\right)\right|  \le C\    x^{2(1+\nu) -  {\rm Re} s}  (x^2- a^2)^{{\rm Re} s/2 -1-\nu } $$

$$\times \left|  {\Gamma(2+\nu)\over \Gamma( (2- s)/ 2) \Gamma( (2(1+\nu)+ s)/ 2)} \right|.$$
These estimates allow to prove the convergence of the integral (11) as an improper one. In fact, writing it as 

$$\lim_{N\to \infty} \int_0^N \varphi(\lambda) \left[ J_\nu(x\lambda) Y_{\nu+1} (a\lambda) -   Y_\nu(x\lambda) J_{\nu+1} (a\lambda)\right] d\lambda,$$
we take  $\varphi$  from the subspace ${\mathcal M}_{1/2,1} ^{-1}(L_c) \subset {\mathcal M}_{0,1}^{-1}(L_c) $ with $c =\{s \in \mathbb{C}:  -1< {\rm Re} s < 0\}$.  So, according to Definition 2 $\varphi$ is given by integral (8) of some function $\Phi(s)$ from the weighted $L_1$-space. Then changing the order of integration by Fubini's theorem for each fixed $N$ and using (14), it becomes 

$$ \int_0^N \varphi(\lambda) \left[ J_\nu(x\lambda) Y_{\nu+1} (a\lambda) -   Y_\nu(x\lambda) J_{\nu+1} (a\lambda)\right] d\lambda=  {1\over 2\pi i}  \int_{\mu- i\infty}^{\mu+i\infty} \Phi(s)  F_\nu(x,s)  ds$$

$$-   {1\over 2\pi i}  \int_{\mu- i\infty}^{\mu+i\infty} \Phi(s)   \int_N^\infty \lambda^{-s} \left[ J_\nu(x\lambda) Y_{\nu+1} (a\lambda) -   Y_\nu(x\lambda) J_{\nu+1} (a\lambda)\right] d\lambda ds,\ x >a.$$ 
We will prove that 

$$\lim_{N\to \infty} \int_{\mu- i\infty}^{\mu+i\infty} \Phi(s)   \int_N^\infty \lambda^{-s} \left[ J_\nu(x\lambda) Y_{\nu+1} (a\lambda) -   Y_\nu(x\lambda) J_{\nu+1} (a\lambda)\right] d\lambda ds=0,\ x > a.\eqno(17)$$ 
To do this, we appeal to the asymptotic behavior of Bessel functions at infinity \cite{erd}, Vol. II to find for fixed $x> a$

$$  J_\nu(x\lambda) Y_{\nu+1} (a\lambda) -   Y_\nu(x\lambda) J_{\nu+1} (a\lambda) = -   {2\over \pi\lambda \sqrt{x a}}
\left[ \cos\left( \lambda(x-a) \right) +  O\left( {1\over \lambda}\right)\right],\  \lambda \to \infty.$$
Substituting this expression into (17) and integrating by parts in the inner integral with respect to $\lambda$ it gives

$$\int_N^\infty \lambda^{-s} \left[ J_\nu(x\lambda) Y_{\nu+1} (a\lambda) -   Y_\nu(x\lambda) J_{\nu+1} (a\lambda)\right] d\lambda =  O\left( (|s| + 1) N^{-\mu -1}\right),$$
and

$$\left| \int_{\mu- i\infty}^{\mu+i\infty} \Phi(s)   \int_N^\infty \lambda^{-s} \left[ J_\nu(x\lambda) Y_{\nu+1} (a\lambda) -   Y_\nu(x\lambda) J_{\nu+1} (a\lambda)\right] d\lambda ds\right|$$

$$ \le C\  N^{-\mu -1}  \int_{\mu- i\infty}^{\mu+i\infty} \left| \Phi(s)\right|  (|s| + 1) |ds| \to 0,\ N \to \infty$$
under assumption $\mu+1 >0.$ Hence we proved the equality 

$$\int_0^\infty  \varphi(\lambda) \left[ J_\nu(x\lambda) Y_{\nu+1} (a\lambda) -   Y_\nu(x\lambda) J_{\nu+1} (a\lambda)\right] d\lambda =   {1\over 2\pi i}  \int_{\mu- i\infty}^{\mu+i\infty} \Phi(s)  F_\nu(x,s)  ds,\eqno(18)$$
where the integral in the right-hand side of (18) converges absolutely. Indeed, from (14), (16) and Stirling's asymptotic formula for the gamma-function at infinity \cite{erd}, Vol. I, we have 

$$\left|F(x,s)\right| \le C\  x^{\mu- {\rm Re} \nu}   e^{\pi |s|/2} |s|^{-\mu},\ x > a.$$
Therefore, 

$$ \int_{\mu- i\infty}^{\mu+i\infty} \left|\Phi(s)  F_\nu(x,s)  ds\right| \le  C\  x^{\mu- {\rm Re} \nu} \int_{\mu- i\infty}^{\mu+i\infty} \left|\Phi(s)\right| e^{\pi |s|/2} |s ds|= C ||\varphi||_{{\mathcal M}_{1/2,1} ^{-1}(L_c)} x^{\mu- {\rm Re} \nu}.$$
Moreover, it tends to zero when $x \to \infty$ when $\mu- {\rm Re} \nu < 0$.   

Let $f \in  {\mathcal M}_{0,1} ^{-1}(L_c) $ with $c =\{s \in \mathbb{C}:   {\rm Re} s =\gamma  > 1/2 \}$.  Returning to integral equation (11) and observing that due to the absolute and uniform convergence of the integral (8) and its derivative with respect to $x \ge x_0 >0$  functions $f,\varphi$ are continuously differentiable on  $[a,\infty)$ and $\mathbb{R}_+$, respectively, we  act  with the differential operator $\left[x^\nu {d\over dx}  x^{-\nu} \right] $ on its both  sides.  Hence, employing  (12), (13), we  obtain

$$\int_0^\infty \lambda \varphi(\lambda) \left[Y_{\nu+1} (x\lambda) J_{\nu+1} (a\lambda)  - J_{\nu+1} (x\lambda) Y_{\nu+1} (a\lambda) \right] d\lambda =  \left[x^\nu {d\over dx}  x^{-\nu} \right] f(x).\eqno(19)$$
It is allowed owing to the uniform convergence by $x \ge a_0 > a$ of the integral with respect to $\lambda$  in (19) if we keep function $\lambda\varphi(\lambda) $ in the same space ${\mathcal M}_{0,1}^{-1}(L_c) $, i.e.  

$$\lambda \varphi(\lambda) =  {1\over 2\pi i}  \int_{\mu- i\infty}^{\mu+i\infty} \Psi(s)  \lambda^{-s} ds,\eqno(20)$$
where $|s| \Psi(s) \in L_1(c),\  c =\{s \in \mathbb{C}:   -1< \mu < 0 \}$.  In fact, recalling the asymptotic behavior of Bessel functions at infinity and integrating by parts, we find 

$$\left|\int_N^\infty \lambda \varphi(\lambda) \left[Y_{\nu+1} (x\lambda) J_{\nu+1} (a\lambda)  - J_{\nu+1} (x\lambda) Y_{\nu+1} (a\lambda) \right] d\lambda\right|$$

$$\le   {1\over \pi^2 \sqrt{x a}} \left|\int_N^\infty  \sin \left( \lambda(x-a) \right) \int_{\mu- i\infty}^{\mu+i\infty} \Psi(s)  \lambda^{-s-1} ds d\lambda\right|   + C N^{-\mu-1}  \int_{\mu- i\infty}^{\mu+i\infty} \left|\Psi(s) ds \right| $$

$$\le  O(N^{-\mu-1} ) +   {1\over \pi^2 (a_0-a) \sqrt{a_0 a}} \left|\int_N^\infty  \cos \left( \lambda(x-a) \right) \int_{\mu- i\infty}^{\mu+i\infty} \Psi(s) (s+1)  \lambda^{-s-2} ds d\lambda\right|$$

$$\le O(N^{-\mu-1} ) +   {N^{-\mu-1} \over \pi^2 (a_0-a) \sqrt{a_0 a}}  \int_{\mu- i\infty}^{\mu+i\infty} |\Psi(s)|  (|s|+1)  | ds| \to 0,\ N \to \infty,$$
where the differentiation under the integral sign is allowed via the absolute and uniform convergence.

Meanwhile, by the same reasons for $x\ge a$, we have (see (8)) 

 $$ \left[x^\nu {d\over dx}  x^{-\nu} \right] f(x)= - {1\over 2\pi i}  \int_{\gamma - i\infty}^{\gamma +i\infty} (s+\nu ) F(s)  x^{-s-1} ds,\ x \ge a,$$
and it tends to zero when $x \to \infty$ via the estimate

$$\left| \left[x^\nu {d\over dx}  x^{-\nu} \right] f(x) \right|\le  
C\  x^{-\gamma -1}  \big\vert\big\vert f  \big\vert\big\vert_{{\mathcal M}_{0,1}^{-1}(L_c)},\  \gamma > {1\over 2}$$
as well as 

$$\left| {d\over dx}  \left[ x^{-\nu}  f(x)\right]  \right|\le  
C\  x^{-\gamma -\nu - 1}  \big\vert\big\vert f  \big\vert\big\vert_{{\mathcal M}_{0,1}^{-1}(L_c)} \to 0,\ x \to \infty.$$
This means that equations (11), (19) are equivalent.  Hence employing Theorem 1, the unique solution of equation (19) has the form

$$  \varphi(\lambda)  =  -   {1\over  J_{\nu+1}^2(a\lambda)+ Y_{\nu+1} ^2(a\lambda )} \int_a^\infty C_{\nu+1}  (\lambda t, \lambda a)\  t^{\nu+1}  {d\over dt}  \left[ t^{-\nu}  f(t) \right] dt,\ \lambda  >0.\eqno(20)$$  
It can be written in a different form with the integration by parts, eliminating outer integrated terms since $t^{1/2} f(t)= o(1),\ t \to \infty$ and taking into account (12), (13).   Hence,

$$ \varphi(\lambda)  =    {1\over  J_{\nu+1}^2(a\lambda)+ Y_{\nu+1} ^2(a\lambda )} \int_a^\infty t  f(t)  \left [  t^{-1-\nu}  {d\over dt}   t^{1+\nu} \right] C_{\nu+1}  (\lambda t, \lambda a) dt$$  

 $$=  {\lambda \over  J_{\nu+1}^2(a\lambda)+ Y_{\nu+1} ^2(a\lambda )} \int_a^\infty t  f(t)  \left[ J_\nu(\lambda t) Y_{\nu+1} (a\lambda) -   Y_\nu(\lambda t) J_{\nu+1} (a\lambda)\right] dt,\ \lambda > 0.$$

We summarize our results by the following

{\bf Theorem 2}. {\it Let  $-1<  {\rm Re } \nu < -1/2$,\  $f \in  {\mathcal M}_{0,1} ^{-1}(L_c) $ with $c =\{s \in \mathbb{C}:   {\rm Re} s =\gamma  > 1/2 \},  \  \varphi \in {\mathcal M}_{1/2,1} ^{-1}(L_c) $ and $\lambda \varphi(\lambda)  \in  {\mathcal M}_{0,1} ^{-1}(L_c) $ with $c =\{s \in \mathbb{C}:  -1< {\rm Re} s <  {\rm Re } \nu \}$.  Then $\varphi$ is the unique solution of the Weber type integral equation $( 11)$, given by the formula}  

$$\varphi(\lambda) =  {\lambda \over  J_{\nu+1}^2(a\lambda)+ Y_{\nu+1} ^2(a\lambda )} \int_a^\infty t  f(t)  \left[ J_\nu(\lambda t) Y_{\nu+1} (a\lambda) -   Y_\nu(\lambda t) J_{\nu+1} (a\lambda)\right] dt,\ \lambda > 0.$$

\bigskip
\centerline{{\bf Acknowledgments}}
\bigskip
The work was partially supported by CMUP [UID/MAT/00144/2013], which is funded by FCT(Portugal) with national (MEC) and European structural funds through the programs FEDER, under the partnership agreement PT2020.   The author thanks Mark Craddock for pointing out the Weber type equation for possible applications  to his attention. 

\bigskip

\bibliographystyle{amsplain}
{}

\end{document}